\documentclass[10pt,twocolumn,letterpaper]{article}

\usepackage{3dv}
\usepackage{times}
\usepackage{epsfig}
\usepackage{graphicx}
\usepackage{amsmath}
\usepackage{amssymb}
\usepackage{bm}
\newcommand*{\bb}[1]{\bm{\mathrm{#1}}}

\usepackage{latexsym}
\usepackage{amsfonts}
\usepackage{color}
\usepackage{comment}
\usepackage[abs]{overpic}
\usepackage{array}
\usepackage{multicol}
\usepackage{multirow}
\usepackage{cite}
\usepackage[boxed]{algorithm2e}
\usepackage{caption}
\usepackage{subcaption}
\usepackage{multicol}
\usepackage{multirow}
\usepackage{array}
\usepackage{dsfont}

\usepackage[pagebackref=true,breaklinks=true,letterpaper=true,colorlinks,bookmarks=false]{hyperref}

\threedvfinalcopy 

\def\threedvPaperID{121} 

\ifthreedvfinal\fi
\begin{document}

\title{Consistent Discretization and Minimization of the $L_1$ Norm on Manifolds}

\author{
 Alex Bronstein \hspace{0.5cm} Yoni Choukroun \cr
     Ron Kimmel \hspace{2.55cm} Matan Sela\cr
     \\
     Computer Science Department \cr
	Technion - Israel Institute of Technology\cr
{\tt\small \{bron,yonic,ron,matansel\}@cs.technion.ac.il}
}

\maketitle

\begin{abstract}
The $L_1$ norm has been tremendously popular in signal and image
 processing in the past two decades due to its sparsity-promoting properties.
More recently, its generalization to non-Euclidean domains has been found
 useful in shape analysis applications.
For example, in conjunction with the minimization of the Dirichlet energy,
 it was shown to produce a compactly supported quasi-harmonic orthonormal basis,
  dubbed as compressed manifold modes \cite{neumann2014}.
The continuous $L_1$ norm on the manifold is often replaced
 by the vector $\ell_1$ norm applied to sampled functions.
We show that such an approach is incorrect in the sense that it does
 not consistently discretize the continuous norm and warn against its sensitivity to the specific sampling.
We propose two alternative discretizations resulting in an
 iteratively-reweighed $\ell_2$ norm.
We demonstrate the proposed strategy on the compressed modes problem,
 which reduces to a sequence of simple eigendecomposition problems not
 requiring non-convex optimization on Stiefel manifolds and producing more stable and accurate results.
\end{abstract}


\section{Introduction}




The $\ell_1$ norm plays a cardinal role in modern digital signal and image processing,
 mainly due to its sparsity-promoting properties and convexity.
Robust PCA \cite{Zou06}, compressed sensing, and inverse problem regularization
 using synthesis and analysis sparse models are just a few examples of
 applications of the $\ell_1$ norm.
The $\ell_1$ norm constitutes a convex surrogate to the combinatorial
 $\ell_0$ norm counting the number of non-zero entries in a vector,
 and powerful theoretical results exist showing the equivalence of such
 convex relaxations of intractable $\ell_0$ minimization problems \cite{Elad2010}.

A limited set of methods involving similar $L_1$-regularization have
 recently appeared for functions defined over discrete surfaces in a variety
 of tasks in geometry processing
  \cite{Dobrev2009,Avron:2010:SRS:1857907.1857911,Deng2013-LocalMod,
  kobbelt2000multiresolution}
  and shape analysis \cite{pokrass2013sparse,cosmo2016consistent,neumann2014}.
The compressed manifold modes (CMM) introduced in \cite{neumann2014} are
 an example of localized smooth truncated basis obtained via sparse
 regularization.
Such bases enjoy most of the properties of the extensively used harmonic
 basis (the orthonormal basis diagnoalizing the Laplacian operator),
 and can constitute an alternative thereof in many geometry processing and analysis tasks.

However, directly copying the $\ell_1$ regularization techniques from
 signal processing hides a potential danger.
Most existing $\ell_1$ regularization models are formulated for problems in
 which continuous signals are sampled at a constant and sufficiently high rate.
In these scenarios, functions can be treated as piecewise constant, and
 the $\ell_1$ norm defined by the sum of absolute values of the
 samples approximates well the continuous $L_1$ norm.
This discretization is no longer valid if the samples are non-uniformly
 distributed or their values have a different meaning.
For example, in computer graphics, shapes are frequently represented as
 discrete triangulated meshes constituting a piecewise-linear approximation of
 the underlying continuous surface.
Numerous methods limit the space of functions on the mesh to be piecewise-linear. For example, in order to solve PDEs numerically.
This approach lies at the core of the finite element method
 (FEM) \cite{zienkiewicz1977finite}.
In such cases, blindly applying the vector $\ell_1$ norm to the
 finite-dimensional vector representing a piecewise-linear function on the
 mesh does not correctly discretize the continuous $L_1$ norm and depends
 on the specific sampling and triangulation.

The present paper addresses this issue.
To that end, we make several contributions.
First, we propose a consistent discretization of the $L_1$ norm evaluated
 as a sum of weighted values, where the weights themselves depend on
 the locations and values of the function samples.
This scheme is compatible with the piecewise-linear representation used in FEM.
Second, for optimization problems involving the proposed norm, we propose
 to translate the objective into a tractable weighted $\ell_2$ norm that
 is minimized by an iterative reweighing scheme.
Finally, we demonstrate experimentally the advantages of the new scheme
 as a general framework for solving $\ell_1$ regularization problems on
 discrete surfaces.
Using the compressed manifold modes problem as a case study, we show that
 it can be formulated as a sequence of eigendecomposition problems,
 avoiding altogether non-convex optimization on Stiefel manifold originally
 used in \cite{neumann2014} and gaining orders of magnitude speedup in runtime.
We show that the resulting bases are robust to different triangulation
 and isometries of a given shape.


\section{Consistent $L_1$ norm discretization}
Consider a continuous surface $\mathcal{M}$ discretized as a triangular mesh
 with the vertices $\mathcal{V} = \{x_i\}_{i=1}^n$ and faces
  $\mathcal{F} = \{t_i\}_{i=1}^m$.
A real-valued function $f:\mathcal{M} \rightarrow \mathds{R}$ defined on the surface is observed at the vertices to form a set of samples $\{ f_i \}_{i=1}^n$ represented as the $n$-dimensional vector $\bb{f}$.
The continuous $\ell_1$-norm of $f$ is defined as
\begin{equation}
\label{eqn:l1_cont}
\begin{aligned}
\| f \|_{L_1} = \int_{\mathcal{M}} \left| f(x) \right| da,
\end{aligned}
\end{equation}
where $da$ denotes the standard area element on $\mathcal{M}$.

There are various ways to define a corresponding discrete norm which approximates the continuous one.
A na\"{i}ve approach is to sum the absolute values of the samples,
$$
\| \bb{f} \|_{\ell_1} = \sum_{i=1}^n |f_i|.
$$
However, this discrete norm does not take the scale of $\mathrm{vol}(\mathcal{M})$ into account, and does not correctly discretize the integral in the overwhelmingly typical case where the vertices of the mesh are non-uniformly distributed.
In what follows, we define two possible alternatives which are directly derived from the zeroth- and first-order approximations of functions over the mesh.

The triangular mesh can be split into $n$ Voronoi cells corresponding to each of the vertices, with the corresponding areas $a_i$.
By approximating the function $f$ on the continuous domain as a piecewise-constant function on the mesh, we assume that the function has the fixed value $f_i$ within each Voronoi cell corresponding to the $i$-th vertex.
We can then straightforwardly define, as in \cite{rustamov2011multiscale}, the area-weighted zeroth order discretization of the $L_1$ norm as
\begin{equation}
\begin{aligned}
\| \bb{f} \|_{L_1^{(0)}} =  \sum_{i=1}^n \left| f_i \right| \cdot a_i.
\end{aligned}
\label{eq:l1-0}
\end{equation}

Alternatively, in the first-order approximation, we assume that the functions are linear within each triangle and thus, piecewise linear over the entire mesh. Hence, the value at a given point $x$ lying in the triangle formed by
 the vertices $(x_i,x_j,x_k)$, is a linear interpolation of the three
  values $(f_i,f_j,f_k)$ at the vertices of the triangle.
The interpolation coefficients are the barycentric coordinates of
 the point $x$.
In other words, $f(x) \approx f_i b_i(x) + f_j b_j(x) + f_k b_k(x)$, where the functions $b_i(x)$ are the piecewise-linear \emph{hat} functions defined as
\begin{equation}
b_i(x)  = \left\{
	\begin{array}{ll}
		1  & : x = x_i \\
		0 &  : x \in \mathcal{M} \setminus \mathcal{N}_1(x_i) \\
		\mathrm{linear} & \mbox{on}\, \mathcal{N}_1(x_i),
	\end{array}
\right.
\end{equation}
where $\mathcal{N}_1(x_i)$ denotes the set of triangles adjacent to $v_i$.

\begin{figure}
    \centering

        \centering
        \includegraphics[scale=0.25]{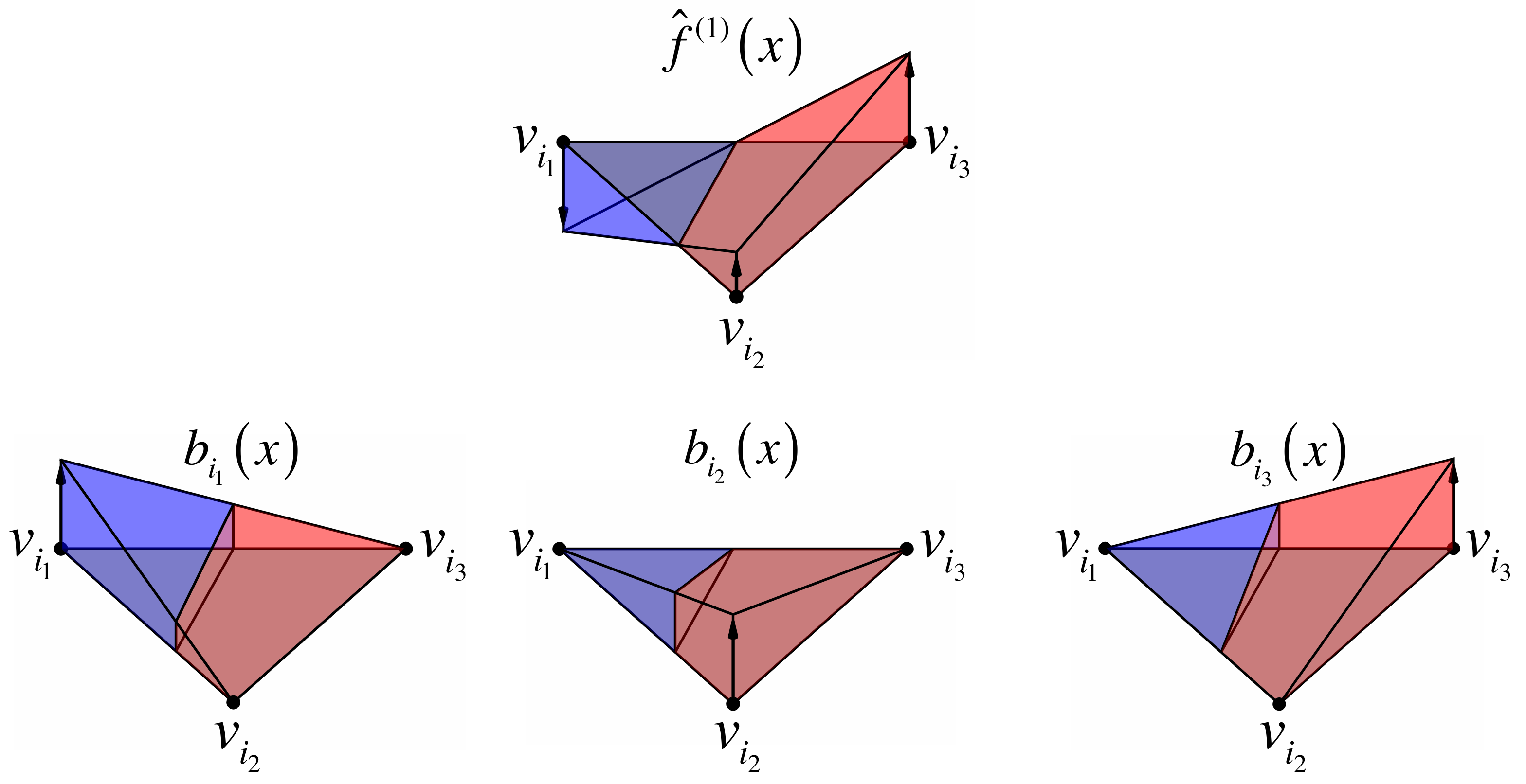}
\caption{Computation of the weight for each vertex in a given triangle. Top: We first find the zero crossing within the triangle. Since the function $\hat{f}^{(1)} (x)$ is linear in each triangle, this is a straight line and splits the triangle into a positive domain (red) and a negative one (blue). (\ref{eq:matan_weight2}). Bottom: The hat basis function of each vertex, colored according to the sign of the function $\hat{f}^{(1)}(x)$.
 In accordance with the integral given in Equation \ref{eq:matan_weight2},  the weight of each vertex in this triangle is computed by subtracting the blue volume from the red one. This weight is summed over all the triangles adjacent to each vertex, for calculating $w_i(\bb{f})$.}
\label{fig:vols}
\end{figure}

Since each $b_i(x)$ is defined over the entire mesh but vanishes outside the $1$-ring of vertex $v_i$, the first order approximation of the function can be written as $ \hat{f}^{(1)}(x) = \sum_{i=1}^n f_i b_i(x) \approx f(x)$.
By plugging the proposed approximation into (\ref{eqn:l1_cont}), a geometric first-order $\ell_1$-norm can be defined as
\begin{eqnarray}\label{eq:matan_weight}
\| \bb{f} \|_{L_1^{(1)}}  &=& \int_{\mathcal{F}} \left| \hat{f}^{(1)}(x) \right| da  =  \int_{\mathcal{F}} \left|  \sum_{i=1}^n f_i b_i(x) \right| da \cr
& = & \sum_{i=1}^n f_i \cdot w_i(\bb{f}),
\label{eq:l1-1}
\end{eqnarray}
where
\begin{eqnarray}\label{eq:matan_weight2}
w_i(\bb{f}) &=& \int_{\mathcal{F}} b_i(x)  \cdot \mathrm{sign} (\hat{f}^{(1)}(x))da \nonumber\\
&=& \sum_{t_j \in \mathcal{N}_1(x_i) } \int_{t_j} b_i(x)  \cdot \mathrm{sign} (\hat{f}^{(1)}(x))da
\end{eqnarray}
and $\mathrm{sign}$ denotes the signum function.
Since, the function is linear in each triangle, the above integrals can be computed simply by calculating volumes of simple polyhedra, as visualized in Figure \ref{fig:vols}.

%

\begin{figure}[t]
\centering
\includegraphics[width=1\columnwidth]{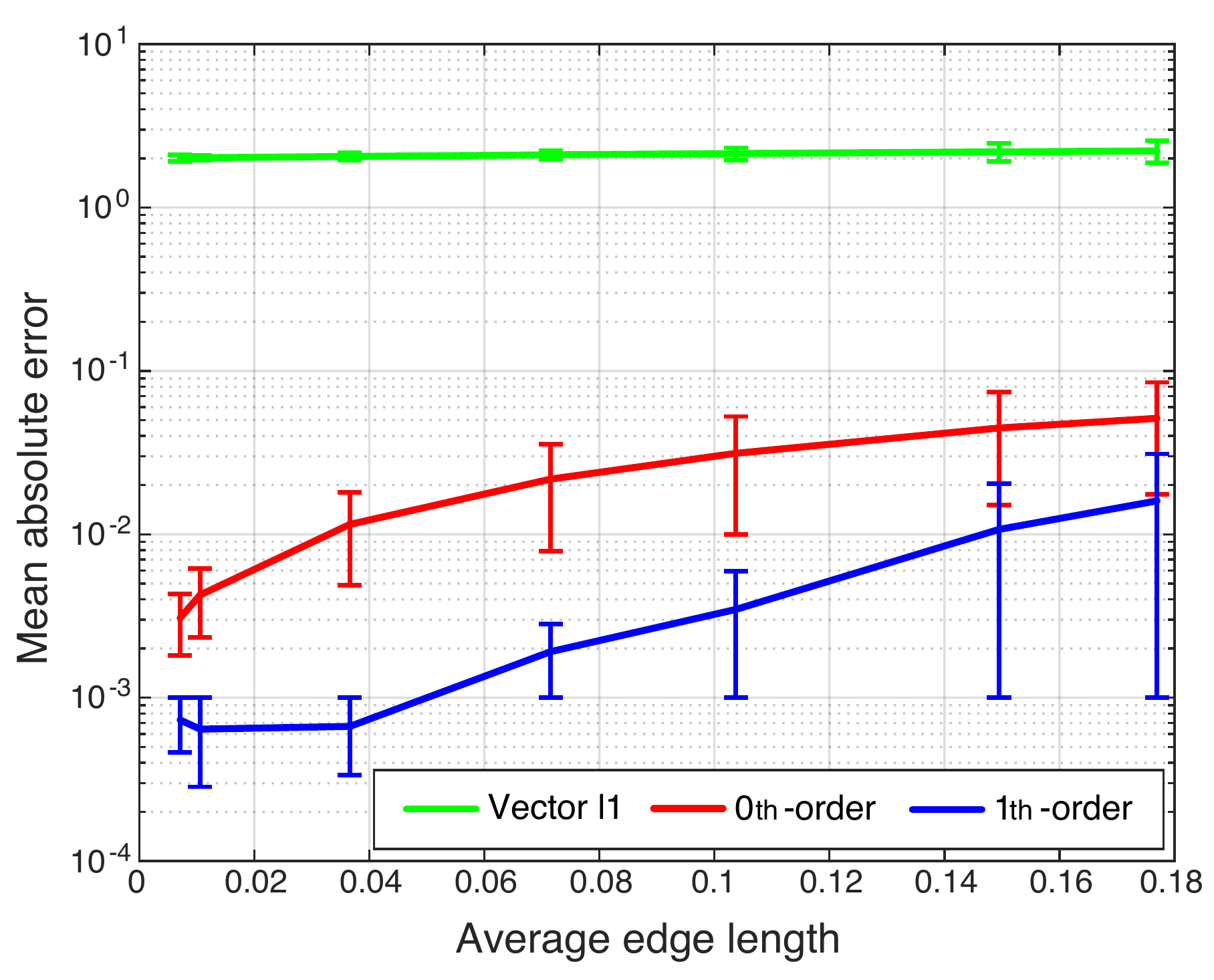}
\caption{Mean absolute approximation error of different $L_1$ norm discretizations plotted as function of mesh resolution (increasing with the decrease of the average edge length). Plotted are the na\"{i}ve vector $\| \cdot \|_{\ell_1}$ norm in green, the zeroth-order approximation $\| \cdot \|_{L_1^{(0)}}$ (\ref{eq:l1-0}) in red, and the first-order approximation $\| \cdot \|_{L_1^{(1)}}$ (\ref{eq:l1-1}) in blue. \label{fig:disc-error}}
\end{figure}

The mean absolute approximation error of different $L_1$ norm discretizations is shown in Figure \ref{fig:disc-error}. As the test functions, we used the first $200$ eigenfunctions of the Laplace-Beltrami operator of a triangular mesh that was remeshed to different resolutions. The area-weighted $\ell_1$ norm of the densely oversampled mesh was used as the reference for error computation, since for sufficiently dense mesh discretization, the difference between the two proposed approximations is negligible.

%
%

\subsection{$L_1$ norm minimization}

One of the bold uses of the $L_1$ norm is its inclusion as a sparsity-promoting penalty or regularization term in optimization problems, giving rise to problems of the form
$$
\min_{\bb{f}}\, E(\bb{f}) + \| \bb{f} \|_1,
$$
where $E(\bb{f})$ is some objective. The lack of smoothness of the norm usually requires the utilization of non-smooth optimization techniques such as proximal \cite{parikh2014proximal} or ADMM  \cite{boyd2011distributed} algorithms.

Here we propose a generic approach to such problems combining the discretization presented above with the well-known iteratively-reweighted least squares (IRLS) method \cite{holland1977robust}. As we show in the following section, such a formulation appears to be beneficial in some problems.

Consider the problem with a discretized $L_1$ term of the form
\begin{equation}
\min_{\bb{f}}\, E(\bb{f}) + \bb{w}(\bb{f})^\mathrm{T} \bb{f},
\label{eq:min-wl1}
\end{equation}
wherein, for the area wieghted zeroth-order norm, $ w_i (\bb{f}) = a_i \cdot sign(f_i)$, and for the first-order one as in \ref{eq:matan_weight2}.
Assuming the weights $\bb{w}(\bb{f})$ are fixed, we can formulate another problem with a weighted $\ell_2$ term of the form
\begin{equation}
\min_{\bb{f}}\, E(\bb{f}) + \bb{f}^\mathrm{T} \bb{C} \bb{f}
\label{eq:min-wl2}
\end{equation}
with the diagonal matrix $\bb{C} = \mathrm{diag}\{ c_1,\dots,c_n\}$ containing the weights.
We would like the two problems to have the same minimizer.
From first order optimality conditions, we require the gradients of $\bb{w}^\mathrm{T} \bb{f}$ and $\bb{f}^\mathrm{T} \bb{C} \bb{f}$ with respect to $\bb{f}$ to vanish at the same point, which yields
\begin{equation}
c_i = \frac{w_i(\bb{f})}{2 f_i}.
\label{eq:w-update}
\end{equation}
The minimization proceeds by solving a sequence of problems of the form (\ref{eq:min-wl2}), each time recalculating the weights according to (\ref{eq:w-update}).

In many cases, the objective $E(\bb{f})$ is a convex quadratic function of the form
\begin{eqnarray}
E(\bb{f}) = \bb{f}^T\bb{Q}\bb{f} + 2\bb{q}^T\bb{f} + c,
\end{eqnarray}
where $\bb{Q}$ is an $n \times n$ symmetric positive semidefinite matrix (often sparse), $\bb{q}$ is an $n$-dimensional vector, and $c$ is a constant.
For guaranteeing a unique solution to problem (\ref{eq:min-wl2}), one must ensure that the matrix $\bb{B} = \bb{Q} + \bb{C}$ is strictly positive definite.
While it is almost everywhere true by construction for the zeroth-order approximation of the $L_1$ norm (except a measure zero set of points where the $f_i$'s vanish), it is not generally so for the first-order approximation, as $w_i(\bb{f})$ and $f_i$ can have opposite signs.

As a remedy, we propose two possible modifications to the matrix $\bb{B}$ arising in the combined objective. The first alternative is to project $\bb{B}$ onto the positive semidefinite cone.
This is performed by computing all the negative eigenvalues $\{\lambda_i < 0\}_{i=1}^{k}$ and the corresponding eigenvectors $\{\bb{\phi}_i\}_{i=1}^k$ of the matrix $\bb{B}$  and subtracting $\sum_{i=1}^k \lambda_i \bb{\phi}_i \bb{\phi}_i^T$ from it.
This comes at the expense of high computational complexity and the risk of ending up with a full matrix.
The second alternative is to modify only the diagonal elements of the matrix $\bb{B}$. According to the Ger\^{s}gorin's circle theorem,  in a diagonally dominant matrix with positive diagonal entrees, defined as a matrix $\bb{B}$ in which each diagonal entry $b_{ii}$ is larger than the sum of absolute off-diagonal entrees in the same row, $\sum_{j\neq i} \left| b_{ij} \right|$, is guaranteed to be positive definite.
Hence, we propose to modify the $i$-th diagonal entry of $\bb{B}$, only for rows in which the diagonal elements are not dominant, by adding the negative gap $\sum_{j\neq i} \left| b_{ij} \right| - b_{ii}$. This computationally efficient modification turns the matrix $\bb{B}$ into a positive definite while maintaining its sparsity.

By changing the matrix $\bb{B}$ in order to turn the problem (\ref{eq:min-wl2}) into a convex one, we slightly modify the original problem. Instead of minimizing the original objective, we minimize a surrogate convex function which is an upper bound of the true objective. However, as we observed in our experiments, since in each iteration we recompute the weights $\bb{w} (\bb{f})$, the sequence of solutions to problem \ref{eq:min-wl2} is monotonously decreasing with respect to its value in the objective of problem \ref{eq:min-wl1}.
\section{Compressed manifold modes}

In what follows, we briefly overview the compressed manifold modes problem used as a case study for the proposed $L_1$ norm discretization.
Ozoli{\c{n}}{\v{s}} et al. \cite{Ozolins2013} proposed a general formalism for sparse solutions to a class of physical problems in Euclidean domains. To that end, they modified the construction of the standard harmonic basis that minimizes the Dirichlet energy among all orthonormal bases by adding an $L_1$ regularization term. The resulting quasi-harmonics were dubbed \emph{compressed modes} of the domain and were shown to be compactly supported \cite{brezis1974solutions,barekat2014support}.
Neumann et al. \cite{neumann2014} extended this construction to manifolds, suggesting the following $L_1$ normalized problem
\begin{equation}\label{eq:CMM}
\begin{aligned}
& \underset{\phi_{i}}{\text{min}}
& & \sum_{i} \int_{\mathcal{M}}\left( \langle \phi_i,  \Delta_{\mathcal{M}} \phi_{i} \rangle_{\mathcal{M}} +\mu |\phi_{i}| \right) da, \ \\
& \text{s.t.}
& & \langle \phi_{i},\phi_{j}\rangle_{\mathcal{M}} = \delta_{ij},
\end{aligned}
\end{equation}
where $\Delta_{\mathcal{M}}$ denotes the Laplace-Beltrami operator and $\langle \cdot, \cdot\rangle_{\mathcal{M}}$ is the intrinsic inner product on $\mathcal{M}$. The non-negative parameter $\mu$ controls the relative importance of smoothness expressed as the Dirichlet energy (first term) and localization expressed as the $L_1$ norm (second term).

Neumann et al. discretized the problem using the na\"{i}ve vector $\ell_1$ norm, obtaining
\begin{equation}\label{eq:CMM-discrete}
\begin{aligned}
& \underset{\bb{\Phi}}{\text{min}}
& & \bb{\Phi}^\mathrm{T} \bb{W} \bb{\Phi} + \mu \| \bb{\Phi} \|_1 \ \\
& \text{s.t.}
& & \bb{\Phi}^\mathrm{T} \bb{A} \bb{\Phi} = \bb{I},
\end{aligned}
\end{equation}
where $\bb{W}$ is the cotangent weight matrix used in the popular Laplacian discretization scheme \cite{pinkall1993computing,meyer2003discrete}.
The non-convex orthogonality constraint combined with the non-smooth objective required the use of non-trivial optimization technique based on ADMM \cite{boyd2011distributed,kovnatsky2015madmm} and proximal operators, guaranteeing no global solution.
The complexity of the compressed modes problem ($\mu>0$) is strikingly higher than the computation of the regular harmonic basis ($\mu=0$) obtained by the simple generalized eigendecomposition $\bb{W}\bb{\Phi} = \bb{A} \bb{\Phi}\bb{\Lambda}$.

\subsection{Iterative reweighting scheme}\label{sec:Iterative_Scheme}

\begin{figure}
    \centering
    \begin{subfigure}[b]{0.25\textwidth}
        \centering
        \includegraphics[scale=0.2,trim={350 40 300 50},clip]{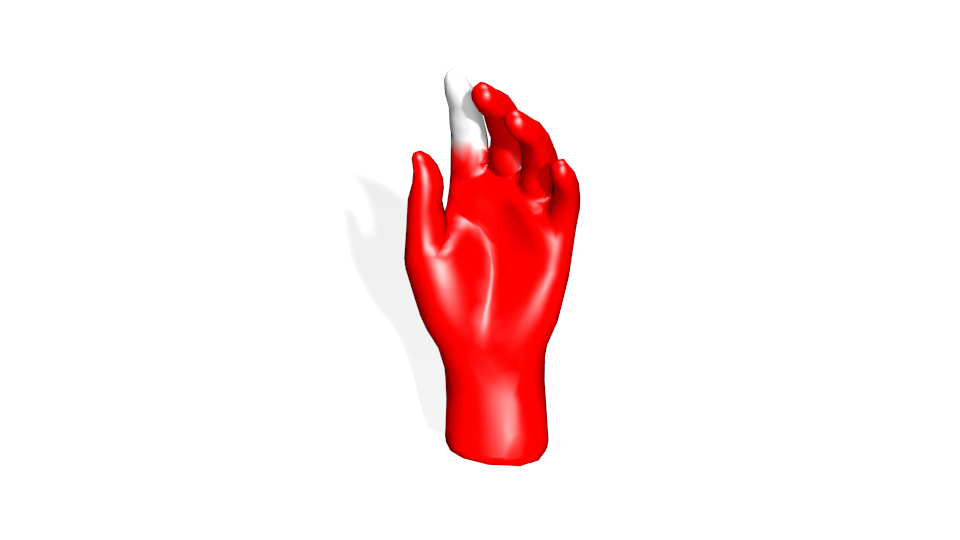}
        \caption{Potential}
    \end{subfigure}%
    ~
    \begin{subfigure}[b]{0.25\textwidth}
        \centering
        \includegraphics[scale=0.2,trim={350 40 300 50},clip]{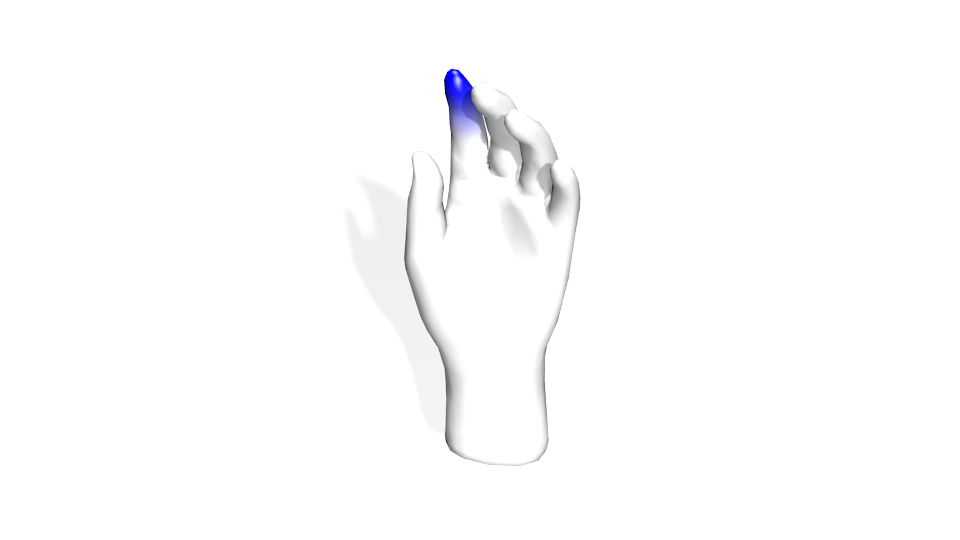}
        \caption{Eigenfunction}
    \end{subfigure}
    \caption{ Potential $v_{i}$ (a), and its corresponding eigenfunction $\phi_{i}$ (b) computed using the proposed framework. Hot and cold colors represent positive and negative values, respectively, while white values represent zero.}\label{fig:potential}
\end{figure}

Using the proposed iteratively-reweighted $L_2$ formulation, we can rewrite the original variational problem (\ref{eq:CMM}) as
\begin{equation}\label{eq:CMM-wl2}
\begin{aligned}
& \underset{\phi_{i}}{\text{min}}
& & \sum_{i} \int_{\mathcal{M}}\langle \phi_i,  \left( \Delta_{\mathcal{M}}  + \mu v_i\right) \phi_{i} \rangle_{\mathcal{M}}  da, \ \\
& \text{s.t.}
& & \langle \phi_{i},\phi_{j}\rangle_{\mathcal{M}} = \delta_{ij},
\end{aligned}
\end{equation}
where $v_i(x)$ can be interpreted as a potential function enforcing diffusion and localizing the support of $\phi_i$ in low-potential areas (Figure \ref{fig:potential}).
Contrary to the original problem (\ref{eq:CMM}), the above problem has a meaningful physical interpretation from quantum mechanics and still looks like operator eigendecomposition.

Using the cotangent discretization of the Laplacian, we formulate the compressed manifold modes problem as the solution of the problem
\begin{equation}\label{eq:CMM-discrete}
\begin{aligned}
& \underset{\bb{\phi}_i}{\text{min}}
& & \bb{\phi}_i^\mathrm{T} \left( \bb{W} \bb{\phi}_i + \mu \bb{A}\bb{V}_{i} \right) \bb{\phi}_i + \beta \sum_{j < i} \| \bb{\phi}_j^\mathrm{T} \bb{A} \bb{\phi}_i  \|_2^2 \\
& \text{s.t.}
& & \bb{\phi}_i^\mathrm{T}\bb{A} \bb{\phi}_i = 1,
\end{aligned}
\end{equation}
where $\beta$ is a sufficiently large constant such that the third term guarantees that the $i$-th mode $\bb{\phi}_i$ is $\bb{A}$-orthogonal to the previously computed modes $\bb{\phi}_j$, $j<i$. Observe that albeit non-convex, the problem has a closed form global solution, that is the smallest generalized eigenvector $\bb{\phi}_i$ of
\begin{equation}
(\bb{W} + \mu \bb{A}\bb{V}_{i} + \beta \bb{Z}_i) \bb{\phi}_i = \lambda_i \bb{A} \bb{\phi}_i
\label{eq:eigendec-irls}
\end{equation}
with
$$
\bb{Z}_i = \bb{A}\left( \sum_{j<i} \bb{\phi}_j \bb{\phi}_j^\mathrm{T} \right)\bb{A}.
$$
When only the few first compressed modes are required, $\bb{Z}_i$ is low rank and finding the smallest  generalized eigenvector can be solved efficiently since the involved matrix is the sum of a sparse and a low-rank matrix.

Several numerical eigendecomposition implementations use the Arnoldi iteration algorithm to extract the eigenvector associated with the
eigenvalue of largest magnitude. The main computationally demanding operation of this method is the multiplication of the matrix we aim at decomposing by a vector. Largest eigenvectors of sparse matrices can therefore be computed very efficiently.
However, since we are seeking the smallest eigenvector, the core operation is the multiplication by the inverse of the matrix we want to decompose. Solved straightforwardly, the iterative solution can be computationally expensive.

For our configuration, let us consider the matrix $\bb{B}$ such that $\bb{B}=\bb{Q}+\bb{U}\bb{U}^\mathrm{T}$ with $\bb{Q}$ and $\bb{U}\bb{U}^\mathrm{T}$ being, respectively, the sparse and the low-rank matrix from (\ref{eq:eigendec-irls}).
Arnoldi's method for the computation of the smallest eigenvector of $\bb{B}$ proceeds by solving at each iteration $\bb{B}\bb{\phi}^{k+1} = \bb{\phi}^{k}$ for the next iterate $\bb{\phi}^{k+1}$ given the current iterate $\bb{\phi}^k$ and normalizing the result.
The Woodbury identity
$$
(\bb{Q}+\bb{U}\bb{U}^\mathrm{T})^{-1} = \bb{Q}^{-1} - \bb{Q}^{-1} \bb{U} ( \bb{I} + \bb{U}^\mathrm{T} \bb{Q}^{-1}\bb{U} )^{-1} \bb{U}^\mathrm{T} \bb{Q}^{-1}
$$
can be used to compute the inverse of the sum of an invertible matrix $\bb{Q}$ and the outer product of two matrices $\bb{U}$ and $\bb{U}^\mathrm{T}$.

At $k$-th iteration, we first compute $\bb{\psi}^k = \bb{Q}^{-1}\bb{\phi}^k$ by solving the sparse system $\bb{Q}\bb{\psi}^k = \bb{\phi}^k$. Next, we compute $\bb{\xi}^k = \bb{Q}^{-1} \bb{U} ( \bb{I} + \bb{U}^\mathrm{T}\bb{Q}^{-1} \bb{U} )^{-1} \bb{U}^\mathrm{T}$ by solving
another sparse system
$$\bb{Q} \bb{\xi}^k = \bb{U} ( \bb{I} + \bb{U}^\mathrm{T}\bb{Q}^{-1} \bb{U} )^{-1} (\bb{U}^\mathrm{T} \bb{\psi}^k).$$
Finally, substituting these two solutions into the Woodbury identity yields
$\bb{\phi}^{k+1} = \bb{B}^{-1} \bb{\phi}^{k} = \bb{\psi}^k - \bb{\xi}^k$.
Since $\bb{Q}$ is sparse, and $\bb{U}$ has only a few columns, the above computations can be carried out efficiently.

\subsection{Experimental evaluation}
\begin{figure}[t]
    \centering
        \includegraphics[width=\linewidth,trim={0 0 0 0},clip]{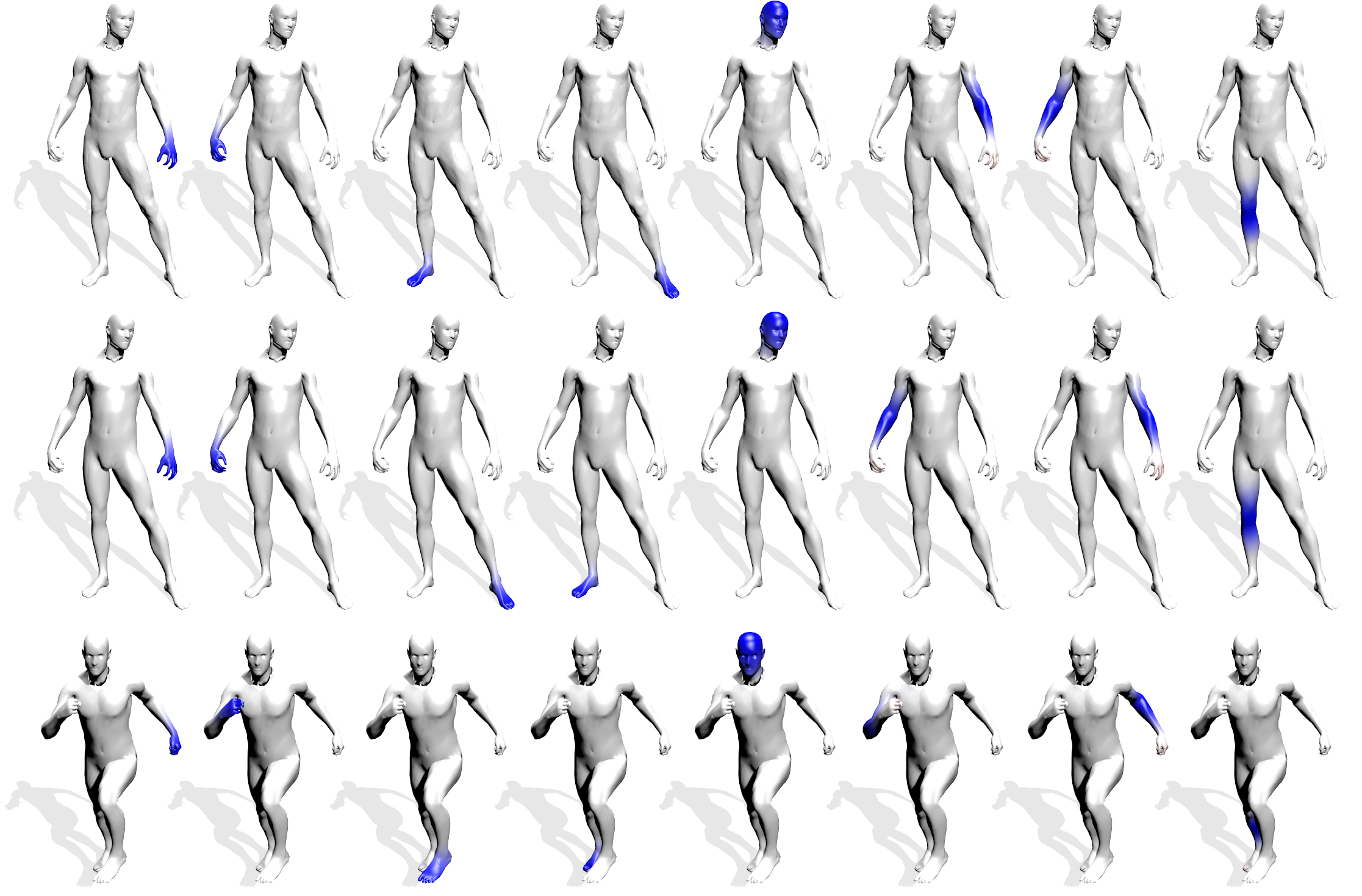}\\
        \noindent\rule{6cm}{0.4pt} \\~\\
          \includegraphics[width=\linewidth,trim={0 0 0 0},clip]{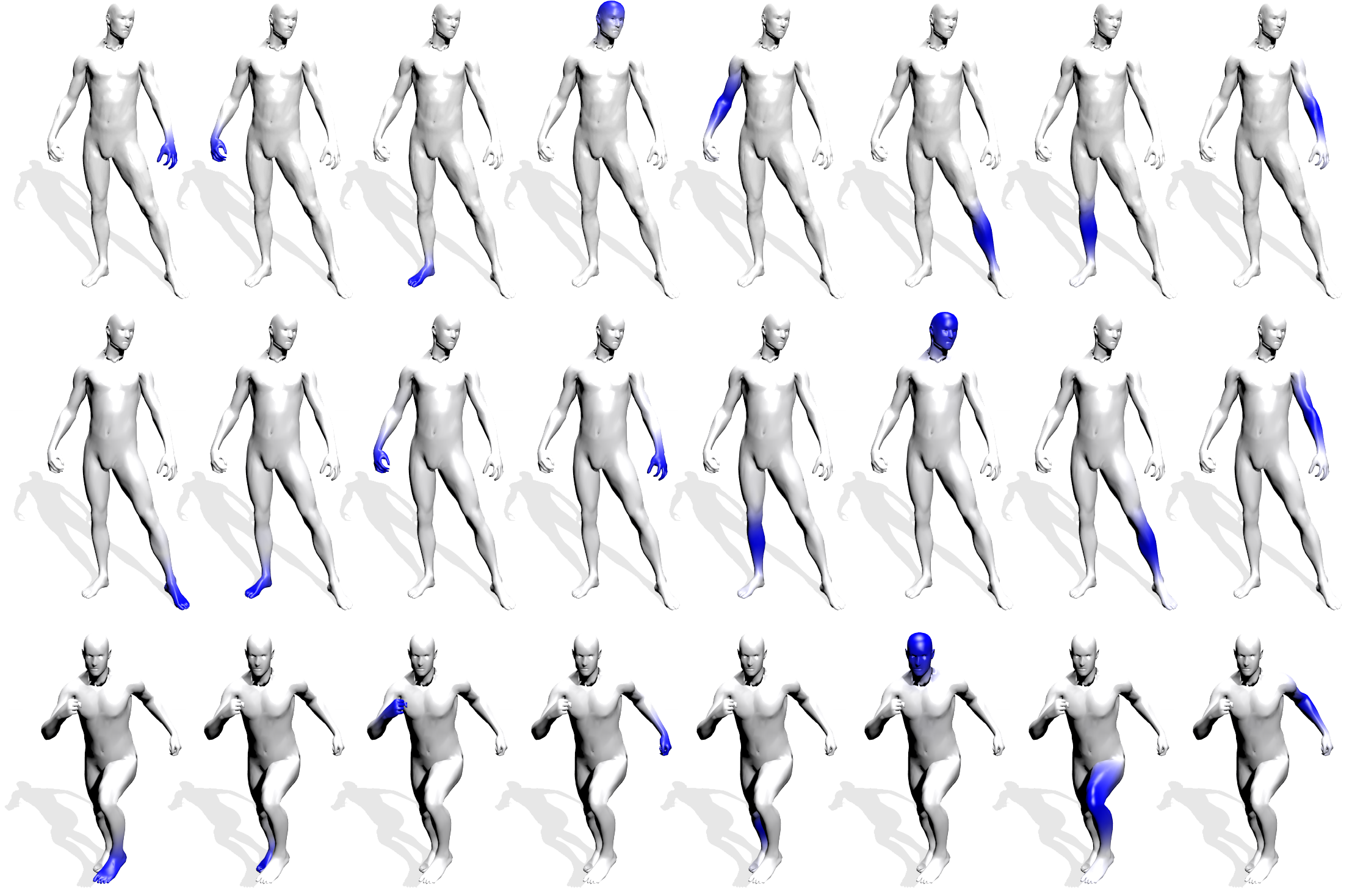}
 \caption{First eight compressed manifold modes computed for different sampling (upsampling by a factor of 5.5) and under nearly isometric deformation of the same mesh using the the proposed IRLS technique (three first rows) and the technique proposed by \cite{neumann2014} (three bottom rows). Note the better stability of the proposed technique to different triangulation.
\label{fig:comp-modes}}
\end{figure}

Figure (\ref{fig:comp-modes}) presents the CMM computed on different surfaces from the TOSCA dataset \cite{bronstein2008numerical} (low resolution) with the method \cite{neumann2014} and the proposed IRLS approach.
The results show spectral decomposition under different sampling, triangulation and deformation. The basis functions obtained by \cite{neumann2014} are sorted according to the cost derived from (\ref{eq:CMM-discrete}). The eigenvectors obtained with the proposed method are naturally sorted by the corresponding eigenvalues. Superior stability under different sampling and nearly isometric deformation of the mesh can be observed.

\begin{figure}[t]
\centering
        \includegraphics[width=\linewidth,trim={0 0 0 0},clip]{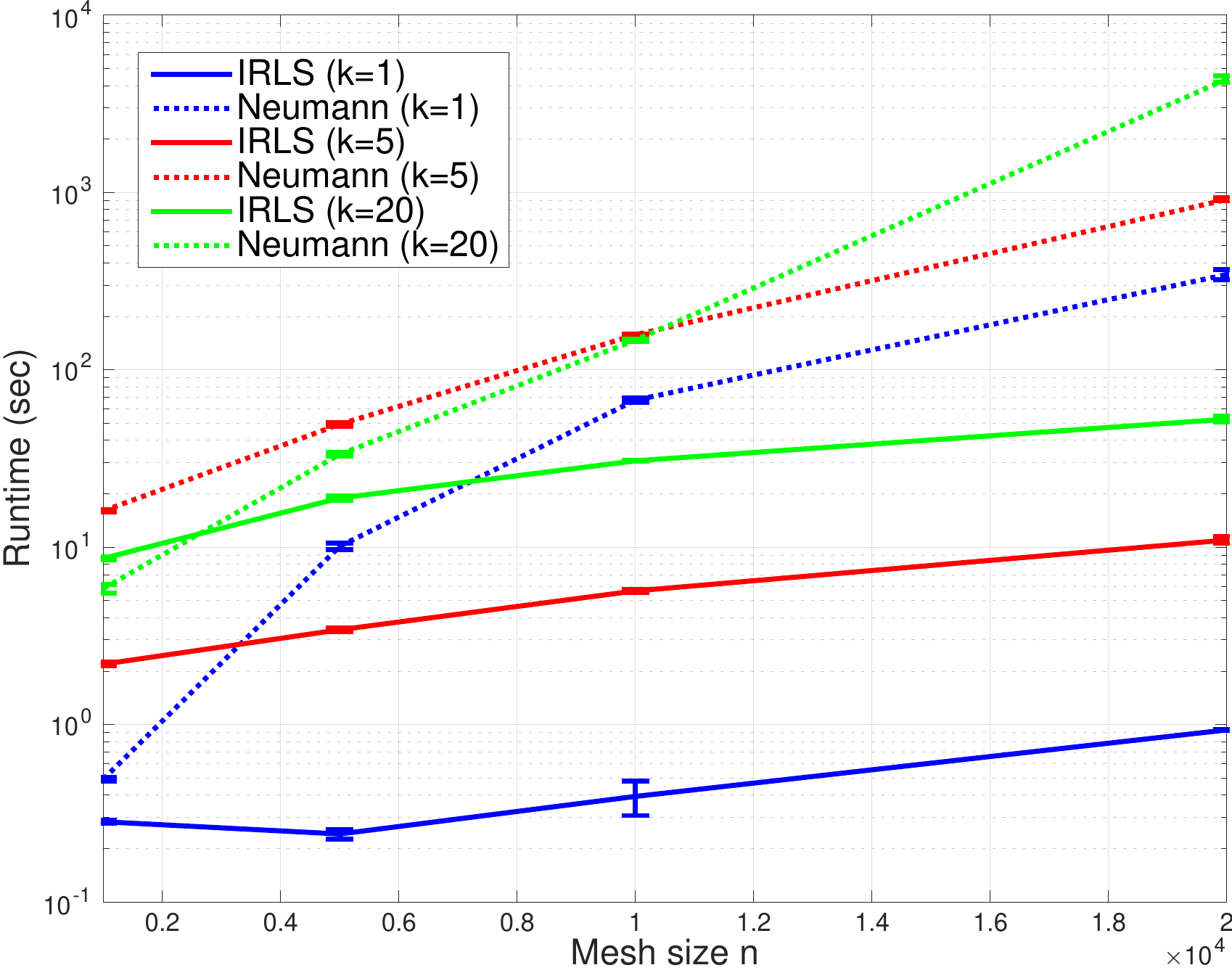}\\
 \caption{Runtimes of \cite{neumann2014} and the proposed framework on meshes of varying size (number of vertices $n$) and number of eigenvectors $k$. Averages and standard deviations are presented over $10$ runs. Same stopping criteria were applied to all methods.
\label{fig:runtimes}}
\end{figure}

For performance comparison we present in Figure (\ref{fig:runtimes}) the runtimes for the different methods using the same sparse parameter $\mu$. The IRLS approach generally requires around 15 iterations to converge, while each iteration is computed efficiently using scheme detailed in Section \ref{sec:Iterative_Scheme}. In other approaches \cite{neumann2014,kovnatsky2015madmm},
high computational complexity make them impractical for dense meshes or when many eigenvectors are required.
The system was implemented in MATLAB and all the experiments
were executed on a 2.5 GHz Intel Core i7 machine
with 16GB RAM. We provided a random initialization
to the reference method \cite{neumann2014} for the computation of the eigenvectors.

\section{Conclusion}
We presented a consistent discretization of the $L_{1}$ norm on manifolds as a geometrically meaningful alternative to the vector $\ell_1$ norm that is frequently employed instead. We also proposed an iteratively-reweighted scheme for the minimization of objectives involving the $L_1$ norm.
As a case study, we considered the recently introduced compressed manifold modes problem and showed that our techniques lead to a significantly more efficient and stable numerical solver.

\section{Acknowledgements}
Y.C, R.K and M.S were supported by Grant agreement No. 267414 of the European Community’s FP7-ERC program. A.B. was supported by the ERC StG Grant 335491 (RAPID). 

{\small
\bibliographystyle{ieee}
\bibliography{references}
}

\end{document}